\documentclass[a4paper]{article}
\usepackage{amsfonts}
\newtheorem{theorem}{Theorem}
\begin{document}
\title{A Nagumo-like uniqueness result for a second order ODE}
\author{
Octavian G. Mustafa\footnote{On leave from the University of Craiova, Faculty of Mathematics and Computer Science, Al{.} I{.} Cuza 13, Craiova, Romania.}\\
\small{University of Vienna, Faculty of Mathematics,}\\
\small{Nordbergstrasse 15, A-1090 Vienna, Austria}\\
\small{e-mail address: octawian@yahoo.com}
}
\date{}
\maketitle

\noindent{\bf Abstract} In this note, we present an extension to second order nonlinear ordinary differential equations (ODEs) of the Nagumo-like uniqueness criterion for first order ODEs established in [A. Constantin, On Nagumo's theorem, Proc. Japan Acad. 86(A) (2010), pp. 41--44].


\noindent{\bf Classification:} 34A12; 34A34

\section{Introduction}
Consider the ordinary differential equation (ODE) below
\begin{eqnarray}
x^{\prime\prime}+f(t,x)=0,\quad t\geq0,\label{main_1}
\end{eqnarray}
where the nonlinearity $f$ is assumed continuous everywhere. In the last few years, many results in the applied sciences were based on sharp uniqueness results for both first order and higher order ODEs dealing with various situations outside the framework of Lipschitzian nonlinearities $f$. See, e{.}g{.}, \cite{kk,Mustafa}.

Our intention in this note is to generalize one of the recent developments in the theory of Nagumo uniqueness criteria \cite{nag,nag2,ath,Mejstrik}, namely the result from \cite{cjapan}. It has been established that for second order ODEs, see (\ref{main_1}), if $\lim\limits_{t\searrow0}f(t,x)=0$ uniformly with respect to $x\in(-1,1)$ and
\begin{eqnarray}
\vert f(t,x_{1})-f(t,x_{2})\vert\leq\frac{2}{t^2}\cdot\vert x_{1}-x_{2}\vert,\quad t>0,\thinspace\vert x_{1}\vert,\vert x_{2}\vert<1,\label{hyp_gen_0}
\end{eqnarray}
the only solution of the equation starting with $x(0)=x^{\prime}(0)=0$ is the trivial one. This has been concluded in \cite{wintner} and further generalized in \cite{c97}. The latter paper presents a Nagumo-like criterion for $n$--th order ODEs when the coefficient $\frac{2}{t^2}$ from (\ref{hyp_gen_0}) is replaced by the Athanassov-like term $\frac{u^{\prime\prime}(t)}{u(t)}$ for some smooth function $u:[0,1]\rightarrow[0,+\infty)$ with $u(0)=0$ and $u(t)$, $u^{\prime}(t)>0$ everywhere in (0,1].

On the other hand, in the case of first order ODEs, both the classical Nagumo restriction
\begin{eqnarray*}
\vert f(t,x_{1})-f(t,x_{2})\vert\leq\frac{1}{t}\cdot\vert x_{1}-x_{2}\vert,\quad t>0,\thinspace\vert x_{1}\vert,\vert x_{2}\vert<1,
\end{eqnarray*}
and its powerful generalization due to Athanassov
\begin{eqnarray*}
\vert f(t,x_{1})-f(t,x_{2})\vert\leq\frac{u^{\prime}(t)}{u(t)}\cdot\vert x_{1}-x_{2}\vert,\quad t>0,\thinspace\vert x_{1}\vert,\vert x_{2}\vert<1,
\end{eqnarray*}
have been extended to
\begin{eqnarray*}
\vert f(t,x)\vert\leq\frac{u^{\prime}(t)}{u(t)}\cdot\omega(\vert x\vert),\quad t\in(0,1],\thinspace\vert x\vert<1,
\end{eqnarray*}
for some continuous and increasing function $\omega:[0,1]\rightarrow[0,+\infty)$ which is null in $0$ and positive everywhere else, and also satisfies the integral inequality
\begin{eqnarray}
\int_{0}^{r}\frac{\omega(s)}{s}\leq r,\quad r\in[0,1].\label{ine_const}
\end{eqnarray}

In the next section, we introduce a variant of the inequality (\ref{ine_const}) such that the equation (\ref{main_1}), where (\ref{hyp_gen_0}) is replaced by
\begin{eqnarray}
\vert f(t,x)\vert\leq\frac{2}{t^2}\cdot\omega(\vert x\vert),\quad t\in(0,1],\thinspace\vert x\vert<1,
\end{eqnarray}
will possess as solution starting from $(0,0)$ only the null solution.

\section{The result and its proof}
\begin{theorem}\label{unu}
Assume that
\begin{eqnarray}
\int_{0}^{r}\frac{\omega(\alpha)}{\alpha^{3/2}}d\alpha\leq 4\sqrt{r},\quad r\in[0,1].\label{main_2}
\end{eqnarray}
Then, the only solution of the equation (\ref{main_1}) starting from $x(0)=x^{\prime}(0)=0$ is the trivial one.
\end{theorem}
Remark that (\ref{main_2}) is satisfied by both the classical example $\omega(r)=2r$, see \cite{wintner,c97}, and the recent example produced by Constantin \cite[p. 43]{cjapan}, namely
\begin{eqnarray*}
\omega(r)=r+\frac{r^2}{2}\sin\frac{1}{r}-\frac{r^2}{3},\quad r\in(0,1],
\end{eqnarray*}
since $\omega(r)\leq r+r^{2}$ and
\begin{eqnarray*}
\int_{0}^{r}\frac{\alpha+\alpha^{2}}{\alpha^{3/2}}d\alpha=2\left(1+\frac{r}{3}\right)\sqrt{r}.
\end{eqnarray*}

{\bf Proof of Theorem \ref{unu}.} Suppose for the sake of contradiction that $x(t)$, $t\in[0,T)$, with $T<1$, is a non-trivial solution of the equation starting from null initial data. Given the local behavior of $f$, that is $\lim\limits_{t\searrow0}f(t,x)=0$ uniformly with respect to $x\in(-1,1)$, we deduce that
\begin{eqnarray*}
\lim\limits_{t\searrow0}x^{\prime\prime}(t)=\lim\limits_{t\searrow0}\frac{x^{\prime}(t)}{t}=\lim\limits_{t\searrow0}\frac{x(t)}{t^{2}/2}=0.
\end{eqnarray*}

There exists the continuous, non-trivial function $y:[0,T]\rightarrow[0,+\infty)$ with the formula
\begin{eqnarray*}
y(t)=\left\{
\begin{array}{ll}
\sup\limits_{s\in(0,t]}\max\left\{\frac{\vert x^{\prime}(s)\vert}{s},\frac{\vert x(s)\vert}{s^{2}/2}\right\},\quad t>0,\\
0,\quad t=0.
\end{array}
\right.
\end{eqnarray*}
According to our assumption, $y(t)>0$ on some subinterval of $(0,T]$.

Further, taking into account the initial datum on the derivative, we have
\begin{eqnarray}
x^{\prime}(t)=\int_{0}^{t}f(s,x(s)) ds,\quad t\in[0,T].\label{main_3}
\end{eqnarray}

We shall apply the inequality regarding $f$ to the integro-differential equation (\ref{main_3}). Thus,
\begin{eqnarray}
\vert x^{\prime}(t)\vert&\leq&\int_{0}^{t}\frac{\omega(\vert x(s)\vert)}{s^2}ds=\int_{0}^{t}\frac{1}{s^2}\cdot\omega\left(\frac{\vert x(s)\vert}{s^{2}/2}\cdot\frac{s^2}{2}\right)ds\nonumber\\
&<&\int_{0}^{t}\frac{1}{s^2}\cdot\omega\left(\frac{\varepsilon s^{2}}{2}\right)ds,\quad t\in(0,T_{\varepsilon}].\label{main_4}
\end{eqnarray}
Here, $T_{\varepsilon}\in(0,T)$ is taken small enough such that $\max\limits_{t\in[0,T_{\varepsilon}]}y(t)=\varepsilon<1$. Remark that
\begin{eqnarray*}
\int_{0}^{t}\frac{1}{s^2}\cdot\omega\left(\frac{\varepsilon s^{2}}{2}\right)ds&=&\frac{\varepsilon}{2}\cdot\int_{0}^{\frac{\varepsilon t^{2}}{2}}\frac{\omega(v)}{v}\cdot\frac{dv}{\varepsilon\sqrt{\frac{2v}{\varepsilon}}}=\frac{\sqrt{\varepsilon}}{2\sqrt{2}}\int_{0}^{\frac{\varepsilon t^{2}}{2}}\frac{\omega(v)}{v\sqrt{v}}dv\\
&\leq&\frac{\sqrt{\varepsilon}}{2\sqrt{2}}\cdot4\sqrt{\frac{\varepsilon t^{2}}{2}}=\varepsilon t
\end{eqnarray*}
by means of the change of variable $\frac{\varepsilon s^{2}}{2}=v$.

Returning to (\ref{main_4}), we get
\begin{eqnarray*}
\vert x^{\prime}(t)\vert<\varepsilon t,\quad t\in(0,T_{\varepsilon}],
\end{eqnarray*}
and respectively (recall that $x(0)=0$)
\begin{eqnarray*}
\vert x(t)\vert\leq\int_{0}^{t}\vert x^{\prime}(s)\vert ds<\frac{\varepsilon t^{2}}{2},\quad t\in(0,T_{\varepsilon}],
\end{eqnarray*}
which lead to $y(t)<\varepsilon$ everywhere in $[0,T_{\varepsilon}]$. This is, obviously, a contradiction given the continuity of the function $y(t)$.

The proof is complete.

\end{document}